\documentclass[12pt,reqno]{amsart}

\usepackage[utf8]{inputenc}
\usepackage[english]{babel}
\usepackage{amsmath,amssymb,amsfonts,amsbsy,amsthm,amscd,latexsym,graphicx}
\usepackage{empheq,textcomp}
\usepackage{mathtools}
\usepackage{relsize}

\theoremstyle{definition}
\newtheorem{theorem}{Theorem} [section]
\newtheorem{corollary}[theorem]{Corollary}
\newtheorem{lemma}[theorem]{Lemma}
\newtheorem{proposition}[theorem]{Proposition}
\newtheorem{definition}[theorem]{Definition}

\newtheorem{remark}[theorem]{Remark}
\newtheorem{example}[theorem]{Example}
\newtheorem{conjecture}[theorem]{Conjecture}
\numberwithin{equation}{section}

\newcommand{\N}{\mathbb{N}}

\newcommand{\Z}{\mathbb{Z}}

\newcommand{\CHI}{\hbox{\raise .4ex \hbox{$\chi$}}} 

\newcommand{\Eq}{\, = \,}

\newcommand{\Ge}{\, \ge \,}

\newcommand{\Le}{\, \le \,}

\newcommand{\qeddef}{{\quad $\diamondsuit$}}

\newcommand{\bigabs}[1]{\bigl|#1\bigr|}
\newcommand{\Bigabs}[1]{\Bigl|#1\Bigr|}

\newcommand{\ip}[2]{\langle#1,#2\rangle}
\newcommand{\bigip}[2]{\bigl\langle #1, \, #2 \bigr\rangle}
\newcommand{\Bigip}[2]{\Bigl\langle #1, \, #2 \Bigr\rangle}

\newcommand{\norm}[1]{\|#1\|}
\newcommand{\bignorm}[1]{\bigl\|#1\bigr\|}
\newcommand{\Bignorm}[1]{\Bigl\|\,#1\,\Bigr\|}

\newcommand{\bigparen}[1]{\bigl(#1\bigr)}
\newcommand{\Bigparen}[1]{\Bigl(#1\Bigr)}
\newcommand{\biggparen}[1]{\biggl(#1\biggr)}
\newcommand{\set}[1]{\{#1\}}
\newcommand{\bigset}[1]{\bigl\{#1\bigr\}}
\newcommand{\Bigset}[1]{\Bigl\{#1\Bigr\}}
\newcommand{\Bigsum}{\mathlarger{\mathlarger{\sum}}}

\newcommand{\clspan}{{\overline{\text{span}}}}

\newcommand{\nf}{{\bigset{\frac{x_n}{\|x_n\|}}}}
\newcommand{\sumli}{\sum_{n=1}^\infty}

\begin{document}

\title{frame-normalizable sequences}

\author{Pu-Ting Yu}

\address{School of Mathematics, Georgia Institute of Technology,
Atlanta, GA 30332, USA}

\email{pyu73@gatech.edu}



\begin{abstract}
Let $H$ be a separable Hilbert space and let $\set{x_n}$ be a sequence in $H$ that does not contain any zero elements. We say that $\set{x_n}$ is a \emph{Bessel-normalizable} or \emph{frame-normalizable} sequence if the normalized sequence $\nf$ is a Bessel sequence or a frame for $H$, respectively. In this paper, several necessary and sufficient conditions for sequences to be frame-normalizable and not frame-normalizable are proved. Perturbation theorems for frame-normalizable sequences are also proved. As applications, we show that the Balazs-Stoeva conjecture 
holds for Bessel-normalizable sequences. Finally, we apply our results to partially answer the open question raised by Aldroubi et al.\ 
as to whether the iterative system $\bigset{\frac{A^n x}{\|A^nx\|}}_{n\geq 0,\, x\in S}$ associated with a normal operator $A\colon H\rightarrow H$ and a countable subset $S$ of $H$, is a frame for $H$. In particular, if $S$ is finite, then we are able to show that $\bigset{\frac{A^n x}{\|A^nx\|}}_{n\geq 0,\, x\in S}$ is not a frame for $H$ whenever $\set{A^nx}_{n\geq 0,\,x\in S}$ is a frame for $H$. 
\end{abstract}

\maketitle

\section{Introduction}
\emph{Frames} are generalizations of bases in a Hilbert space $H$ in the sense that every element of $H$ can be expressed as an ``infinite linear combination" of frame elements with unconditional convergence in norm. The concept of frames was introduced by Duffin and Schaefer \cite{DS52} as part of their research on non-harmonic Fourier series. Frames were popularized by Daubechies, Grossmann, and Meyer  \cite{DGM86} in their study of Gabor systems and wavelets. Since then, there has been a plethora of research in both applied and theoretical mathematics, as well as science and engineering, in this field. For relatively recent textbook recountings, we refer to \cite{Chr16}, \cite{Hei11}, and \cite{You01}.

We say that a countable sequence $\set{x_n}$ (implicitly indexed by the set of natural numbers $\N$)
in a separable Hilbert space $H$ is a frame if there exist positive constants $A \le B,$ called \emph{frame bounds},
such that 
\begin{equation}
\label{frame_ineq}
A\, \norm{x}^2
\Le \sum_{n=1}^\infty |\ip{x}{x_n}|^2
\Le B\, \norm{x}^2,
\qquad\text{for all } x \in H.
\end{equation}
A frame is tight if we can take $A=B$, and Parseval if we can take $A=B=1$.
A sequence $\set{x_n}$ is a \emph{Bessel sequence} in $H$ if at least the upper inequality of equation (\ref{frame_ineq}) is satisfied.
We say that a sequence satisfies a lower frame condition or admits a lower frame bound if at least the lower inequality of equation (\ref{frame_ineq}) holds.

Unit-norm tight frames (UNTFs) have proved useful in many applications, especially in the finite-dimensional setting.
Finite unit-norm tight frames (FUNTFs) avoid problems that arise when applying finitization methods to infinite frames. As a result,
FUNTFs have been widely applied to coding theory \cite{HS03}, sparse representations \cite{HHT19}, and many other areas. 
It is known that the construction of UNTFs is difficult even in the finite-dimensional setting (\cite{BC10}, \cite{CK03}). Benedetto and Fickus characterized FUNTFs via a \emph{frame potential} \cite{BF03}. 
In infinite-dimensional Hilbert spaces, Kutyniok et al.\ studied frames which can become a Parseval frame after being rescaled by a sequence of scalars \cite{KOPT13}.

 Disregarding the tightness requirement, is there any quick way to construct a unit-norm frame for $H$? One way to obtain a unit-norm sequence is to normalize a sequences by dividing every element in the sequence by its length,  but then we must see if we can still obtain a frame in this manner.  Our interest here is to find characterizations of sequences which become frames or Bessel sequences after being normalized. 

Another motivation is the following two open questions. Let $A\colon H\rightarrow H$ be a bounded linear operator and let $S$ be an arbitrary countable subset of $H$. In \cite{ACMCP16}, Aldroubi et al.\ proved that the iterative system $\bigset{\frac{A^n x}{\norm{A^n x}}}_{n\geq 0, \,x\in S}$ is never a frame for $H$ if $A$ is a self-adjoint operator. The question whether $\bigset{\frac{A^n x}{\norm{A^n x}}}_{n\geq 0, \,x\in S}$ can be a frame for $H$ when $A$ is a normal operator was asked in the same paper. When $S$ is a finite subset of $H$, we are able to answer this question negatively if $\bigset{A^n x}_{n\geq 0, \,x\in S}$ is a frame for $H$. Several sufficient conditions for which $\bigset{\frac{A^n x}{\norm{A^n x}}}_{n\geq 0, \,x\in S}$ cannot be a frame will also be proved. 

The second open question is the Balazs-Stoeva conjecture from \cite{BS11}. 
 Let $\mathcal{M}=(m_n)$ be a sequence of scalars, and let $\mathcal{X}=\set{x_n}$ and $\mathcal{Y}=\set{y_n}$ be sequences of elements in $H$. \emph{Multipliers} are operators on $H$ that have the form $\mathcal{M}_{\mathcal{X},\mathcal{Y}}\,(x)= \sum m_n\,\ip{x}{y_n}\,x_n$.
 The Balazs-Stoeva conjecture is that $\sum m_n\,\ip{x}{y_n}\,x_n$ converges unconditionally for all $x\in H$ if and only if 
 there exist sequences of scalars $(c_n)$ and $(d_n)$ with $c_n\overline{d_n}=m_n$ such that $\set{c_nx_n}$ and $\set{d_ny_n}$ are Bessel sequences. Using unconditional convergence of frame expansion as a criterion, the Balazs-Stoeva conjecture gives a potential classification of ``good" and ``bad" alternative duals associated with a frame \cite{HY22}. As a corollary of Orlicz's Theorem, we will see that the Balazs-Stoeva conjecture holds for all sequences which become Bessel sequences after being normalized.
 
 This paper is organized as follows. In Section \ref{prelim_sec}, the definitions, notations and some background results needed for this paper are stated. Section \ref{Bess_frame_norm_sequences} presents properties and characterizations of sequences which become Bessel sequences or frames for $H$.
 Perturbation theorems are studied in Section \ref{pert_norm_frames}. Finally, we apply our results to Balazs-Stoeva conjecture and study the frame property of normalized iterative systems $\bigset{\frac{A^n x}{\norm{A^n x}}}_{n\geq 0, \,x\in S}$ in Section \ref{BS_conjecture}.

\section{Preliminaries} \label{prelim_sec}
We state without proof some results that are required for this paper in this section. We refer to \cite{Chr16}, \cite{Hei11}, and \cite{You01} for more details on frames.

Throughout this paper, $H$ will denote a separable Hilbert space equipped with inner product $\ip{\cdot}{\cdot}$.  $\ell^2$ is the space of square-summable sequences of scalars indexed by  the natural number $\mathbb{N}$. The notation $\{x_n\}$ or $(c_n)$ will denote a countable sequence indexed by $\mathbb{N}$. If $S$ is a closed subspace of $H,$ then $P_S$ will denote the
orthogonal projection of $H$ onto $S.$ Without loss of generality, we assume that all sequences $\set{x_n}$ in our discussion do not contain any zero elements.

We first formalize the notion of Bessel-normalizable and  frame-normalizable sequences.
\begin{definition}
	Let $\set{x_n}$ be a sequence in $H$ which does not contain any zero elements. 
\begin{enumerate}
\setlength\itemsep{0.3em}
\item [\textup{(a)}] We say that $\set{x_n}$ is \emph{Bessel-normalizable} if $\nf$ is a Bessel sequence. If $\set{x_n}$ is a Bessel sequence, then we say that $\set{x_n}$ is a \emph{normalizable Bessel sequence} if it is Bessel-normalizable.
\item [\textup{(b)}] We say that $\set{x_n}$ is \emph{lower-normalizable} if $\nf$ satisfies a lower frame condition, i.e., if there exists a constant $A>0$ such that $\sum\,\bigabs{\bigip{x}{\frac{x_n}{\|x_n\|}}}^2\geq A\,\|x\|^2$  for all $x\in H$.
\item [\textup{(c)}] We say that $\set{x_n}$ is \emph{frame-normalizable} if $\nf$ is a frame. If $\set{x_n}$ is a frame, then we say that $\set{x_n}$ is a \emph{normalizable frame} if it is frame-normalizable. \qeddef
\end{enumerate}
\end{definition}

A sequence $\set{x_n}\subseteq H$ is said to be minimal if it does not contain any redundant elements in the following sense. 

\begin{definition}
	A sequence $\set{x_n}\subseteq H$ is \emph{minimal} if $x_m \notin \clspan{\set{x_n}_{n\neq m}}$ for every $m$. 
	Further, $\set{x_n}$ is \emph{exact} if it is both minimal and complete (here completeness means that the closed span is $H$).\qeddef
\end{definition}
Minimality of a sequence $\set{x_n}\subseteq H$ is equivalent to the existence of a biorthogonal sequence $\set{a_n}\subseteq H$ such that $\ip{x_n}{a_k}=\delta_{nk}$; for example, see \cite[Lem. 5.4]{Hei11}.

We will need the following two bounded linear operators induced by a Bessel sequence.
\begin{definition} If $\set{x_n}$ is a Bessel sequence in $H$, then we define the following associated bounded linear operators.
\begin{enumerate}
	\setlength\itemsep{0.3em}
\item [\textup{(a)}] The analysis operator $C\colon H\rightarrow \ell^2$ is $C(x)=(\ip{x}{x_n})$ for $x\in H$.
\item [\textup{(b)}] The frame operator  $S\colon H \rightarrow H$ is $S(x)= \sum \ip{x}{x_n}\,x_n $ for $x\in H$.\qeddef
\end{enumerate}
\end{definition}
If $\set{x_n}$ is a frame, then the range of the associated analysis operator is a closed subspace of $\ell^2$ and $S$ is a bounded linear operator that is positive and invertible.




One way to transform a frame into a Parseval frame is to apply the inverse of the square root of the frame operator; see \cite[Cor. 8.28]{Hei11}.
\begin{proposition}
	\label{trans_frame_into_pars}
	If $\set{x_n}$ is a frame for $H$ with frame operator $S$, then $\set{S^{-1/2}x_n}$ is a Parseval frame for $H$. 
	Consequently, $\|S^{-1/2}x_n\|\leq 1$ for all $n$. \qeddef
	\end{proposition}

We will need an elegant inequality for Parseval frames proved by Balan et al. in \cite[Prop. 4.1]{BGEK07}.
\begin{proposition}
	\label{iden_parseval}
	If $\set{x_n}$ is a Parseval frame for $H$, then for any subset $J\subseteq \mathbb{N}$ and every $x\in H$ we have 
	$$\sum_{n\in J} |\ip{x}{x_n}|^2\,+\,\Bignorm{\sum_{n\in J^c}\ip{x}{x_n}{x_n}}^2\Ge \frac{3}{4} \|x\|^2 $$
	Further, we have equality for all $x\in H$ if and only if $\sum_{n\in J} \,\ip{x}{x_n}\,x_n=\frac{1}{2} x$ for every $x\in H$. \qeddef
	\end{proposition}
 
 In this paper, we refer to a bounded linear invertible operator as a \emph{topological isomorphism}. Let $H$ and $K$ be two Hilbert spaces. Then sequences $\set{x_n}\subseteq H$ and $\set{y_n}\subseteq K$ are said to be \emph{topologically isomorphic} if there is a topological isomorphism $T\colon H\rightarrow K$ such that $T(x_n)=y_n$ for all $n$. 
 
The next result states that frames are topologically isomorphic to the orthogonal projection of the standard basis for $\ell^2$ onto some closed subspace; see \cite[Cor. 8.33]{Hei11} or \cite[Thm. 2.1]{Hol94}.
\begin{theorem}
	\label{fr_equi_orth_basis_el2}
	Assume that $\set{x_n}$ is a sequence of elements in $H$. Then $\set{x_n}$ is a frame if and only if there exists a closed subspace $S$ of $\ell^2$ and a topological isomorphism $T\colon S\rightarrow H$ such that $x_n=TP_S \delta_n$ for every $n$, where $\set{\delta_n}$ is the standard basis for $\ell^2$ and $P_S$ denotes the orthogonal projection of $\ell^2$ onto $S$.
	
		Further, in case that these statements hold we can take the closed subspace $S$ to be $S= \text{range}(C)$, where $C$ is the analysis operator for $\set{x_n}$. \qeddef
\end{theorem}
\section{Bessel-normalizable and frame-normalizable sequences}
\label{Bess_frame_norm_sequences}
\subsection{General Results}
We start this section by giving some examples. 

Frames are bounded above in norm due to the upper frame inequality of equation (\ref{frame_ineq}), but they are not necessarily bounded below in norm. A frame $\set{x_n}$ is said to be norm-bounded below if $\inf \|x_n\|>0$. 

\begin{example}
	\label{nbb_normalizable}
By a straightforward calculation, every frame that is norm-bounded below is frame-normalizable.
\end{example}
There exist frame-normalizable sequences that are not norm-bounded below.

\begin{example}
	Let $\set{e_n}$ be an orthonormal basis for $H$. If $\set{x_n}=\set{e_n}\cup \bigset{\frac{1}{n}e_n}$, then $\set{x_n}$ is a frame with frame bounds $1$ and $2$. After normalization, we have that $\nf=\set{e_n}\cup \set{e_n}$, which is a tight frame with frame bound $2$.
	\qeddef
\end{example}

The following proposition characterizes all Bessel-normalizable sequences in finite-dimensional Hilbert spaces.
\begin{proposition}
	\label{fini_dim_char}
	Let $H$ be a finite-dimensional Hilbert space and let $\set{x_n}$ be a sequence in $H$. Then $\set{x_n}$ is Bessel-normalizable if and only if it is a finite sequence.
\end{proposition}
\begin{proof}
		($\Rightarrow$) Suppose to the contrary that $\set{x_n}$ is not a finite sequence.
		Let $S=\clspan{\set{x_n}}$ and  let $\set{e_n}_{n=1}^N$ be an orthonormal basis for $S$. Then for each $n$, we have $\frac{x_n}{\|x_n\|}=\sum _{i=1}^N c_{i}^ne_i$ for some sequence of scalars $(c^n_i)_i$. Since $\sum _{i=1}^N |c_{i}^n|^2=1$ for every $n$, 
		there exists an index $1\leq j_0\leq N$ such that $|c^n_{j_0}|^2 \geq \frac{1}{N}$ for infinitely many $n$. Thus, $$\sumli\,\Bigabs{\Bigip{e_{j_0}}{\frac{x_n}{\|x_n\|}}}^2 \geq \sumli\, \frac{1}{N} =\infty.$$
		Therefore, $\nf$ is not a Bessel sequence, which is a contradiction. 
	\medskip
	
	\medskip
	($\Leftarrow$) Follows from the Cauchy–Bunyakovsky–Schwarz Inequality.
	\end{proof}

We say that a sequence $\set{x_n}$ in $H$ admits an \emph{orthogonal decomposition} if there exist (finitely or countably many) subsequences $\set{x_n^1},\set{x_n^2},\dots$, such that $$\set{x_n}=\set{x_n^1}\,\cup\, \set{x_n^2}\, \cup \,\cdots$$ and the closed spans of the $\set{x_n^i}_n$'s are mutually orthogonal, i.e., $\clspan\set{x_n^i}_n\perp \clspan\set{x_n^j}_n $ for $i\neq j$.

A sufficient condition for a sequence to be Bessel-normalizable is that it admits a certain type of orthogonal decomposition.

\begin{theorem}
	\label{norm_Bes_orth_de}
	Let $\set{x_n}$ be a sequence in $H$. Assume that $\set{x_n}$ admits an orthogonal decomposition $\set{x_n}=\set{x_n^1} \,\cup \,\set{x_n^2} \cup \cdots$ and assume further that $$\sup_i \bigparen{\text{dim}( \clspan\set{x_n^i}_n)}<\infty.$$ Then $\set{x_n}$ is Bessel-normalizable if and only if $\sup_i |\set{x_n^i}_n|<\infty$, 
	where  $|\set{x_n^i}_n|$ denotes the cardinality of $\set{x_n^i}_n$. 
\end{theorem}
\begin{proof}
	($\Rightarrow$) The idea is similar to the proof of Proposition \ref{fini_dim_char}. Let $L=\sup_i\bigparen{\text{dim}( \clspan\set{x_n^i}_n)}$. Suppose that $\sup_i |\set{x_n^i}_n|=\infty$. Then for any $K\in \mathbb{N}$, there exists some $i_K\in\mathbb{N}$ such that $|\set{x_n^{i_K}}_n|>L^2K$. Let $\set{e_1,e_2,...,e_\ell}$ be an orthonormal basis for $\clspan{\set{x_n^{i_K}}_n}$, and note that $\ell \leq L$. Then $\frac{x_{n}^{i_K}}{\|x_{n}^{i_K}\|}=\sum_{i=1}^\ell c_i^ne_i$ for some scalars $c_i^n$. Since $\sum_{i=1}^\ell |c_i^n|^2 =1$, at least one of $|c_i^n|^2$ is greater than $\frac{1}{\ell}$. Furthermore, there exists some $1\leq j_0\leq \ell$ such that $|c^n_{j_0}|^2\geq \frac{1}{\ell}$ for at least $LK$ distinct values of $n$. Consequently, $$ \Bigsum_{n=1}^{|\set{x_n^{i_K}}_n|}\, \Bigabs{\Bigip{e_{j_0}}{\,\frac{x_{n}^{i_K}}{\norm{x_{n}^{i_K}}}}}^2\Ge \frac{LK}{\ell}\Ge K\Eq K\|e_{j_0}\|^2.$$ Since $K$ is arbitrary, it follows that $\nf$ is not a Bessel sequence.\newline

	($\Leftarrow$) Assume that $\sup_i |\set{x_n^i}_n|<\infty$. 
	Let $S_i=\clspan\set{x_n^i}_n$. Then
	\begin{align*}
	\label{eq2}
	\begin{split}
	\sum_{\|x_n\|< \delta}  \Bigabs{\Bigip{x}{\frac{x_n}{\|x_n\|}}}^2 &\Eq \sum_{i\in \mathbb{N}}\,\sum_{n=1}^{|\set{x_n^i}_n|} \,\Bigabs{\Bigip{P_{S_i}x}{\frac{x_n^i}{\|x_n^i\|}}}^2 \allowdisplaybreaks \\
	&\Le \sup_i |\set{x_n^i}_n|~\sum_{i\in \mathbb{N}} \,\|P_{S_i}x\|^2  \allowdisplaybreaks \\
	&\Le\sup_i |\set{x_n^i}_n| ~\|x\|^2.
	\end{split}
	\end{align*}
	Thus, $\set{x_n}$ is Bessel-normalizable.
\end{proof}

\begin{remark}
	The assumption that $\sup_i \bigparen{\text{dim}( \clspan\set{x_n^i})}<\infty$ is not required for the ``if" direction in Theorem \ref{norm_Bes_orth_de}, but it is necessary for the ``only if" direction (see Example \ref{non_nor_frames}).  Moreover, if $\set{x_n}$ is a frame, then the proof of Theorem \ref{norm_Bes_orth_de} shows that  $\set{x_n}$ is a normalizable frame if and only if $\sup_i |\set{x_n^i}|<\infty$ by Proposition \ref{equiv_nor_Be_fr}. \qeddef
\end{remark}
%

We recall that a sequence is a Riesz basis for $H$ if it is topologically isomorphic to an orthonormal basis for $H$. A sequence $\set{x_n}$ is a Schauder basis for $H$ if for any $x\in H$ there exists a unique sequence of scalars $(c_n)$ such that $x=\sum c_nx_n.$ 
It was conjectured by Feichtinger that every frame that is norm-bounded below is a finite union of Riesz sequences (i.e., Riesz bases for their closed spans). Later, Casazza et al. showed that this conjecture is equivalent to the conjecture that every Bessel sequence that is norm-bounded below is a union of finitely many Riesz sequences \cite{CCLV05}. Following that, Casazza and Tremain showed that the Feichtinger Conjecture is equivalent to some long-standing intractable problems including the Kadison-Singer Conjecture \cite{CT06}. The Kadison-Singer Conjecture was finally confirmed in \cite{MSS15} by Marcus, Spielman, and Srivastava by their proof of another equivalent formulation. Based on Feichtinger's Conjecture (now Feichtinger's Theorem), we obtain the following decomposition of Bessel-normalizable sequences. 

\begin{proposition}
	\label{nor_Bess_Schau_decom}
	Every Bessel-normalizable sequence is a finite union of Schauder sequences (i.e., Schauder bases for their closed spans). 
\end{proposition}
\begin{proof}
	By Feichtinger's Theorem, there exists some $N\in \N$ such that $\nf$ can be decomposed into $\nf=\bigset{\frac{x_n^1}{\|x_n^1\|}}\cup\bigset{\frac{x_n^2}{\|x_n^2\|}}\cup \cdots\cup \bigset{\frac{x_n^N}{\|x_n^N\|}}$ and $\bigset{\frac{x_n^i}{\|x_n^i\|}}$ is a Riesz sequence for each $i$.
	Fix $1\leq i\leq N$ and let $\set{\widetilde{x}^{\,i}_n}$ be the biorthogonal sequence of $ \bigset{\frac{x_n^i}{\|x_n^i\|}}.$ Then every $x\in \clspan{\set{x_n^i}}$ can be expressed uniquely as $$x\Eq\sum\, \ip{x}{\widetilde{x}^{\,i}_n}\,\frac{x_n^i}{\|x_n^i\|}\Eq\sum\, \Bigip{x}{\frac{\widetilde{x}^{\,i}_n}{\|x_n^i\|}} \,x_n^i.$$
	Thus, $\set{x_n^i}$ is a Schauder sequence. Hence every Bessel-normalizable sequence is a finite union of Schauder sequences.
\end{proof}
\begin{remark}
	The converse does not hold in general. Let $H=L^2(\mathbb{T})$ and consider the function $\phi(x)= |x-\frac{1}{2}|^{-1/4}$. It is known the system of weighted exponentials $\mathcal{E}(\phi) = \set{\phi(x) e^{2\pi inx}}_{n\in \Z}$ is a Schauder basis (but not a Riesz basis) for $L^2(\mathbb{T})$ with respect the ordering $\mathbb{Z}=\set{0,1,-1,2,-2,\dots}$ (see \cite[pp.\ 351]{IS70}). Also, it is known that a system of weighted exponentials $\mathcal{E}(\psi)$ is a Bessel sequence if and only if $\psi \in L^\infty (\mathbb{T})$ (for more characterizations of weighted exponential systems, see \cite[Thm. 10.10]{Hei11}). Consequently, $\mathcal{E}(\phi)$ is not a Bessel sequence. Therefore, since $\norm{\phi(x) e^{2\pi inx}}=\norm{\phi}$ for every $n$, $\mathcal{E}(\phi)$ is not Bessel-normalizable.
	
	 Moreover, Schauder bases that are Bessel sequences are not necessarily Bessel-normalizable. Let $(c_n)\in \ell^2$ with $c_n> 0$ for all $n$. Then $\set{c_n\phi(x) e^{2\pi inx}}$ is a Schauder basis (with respect to the same ordering above) that is a Bessel sequence. However, $\set{c_n\phi(x) e^{2\pi inx}}$ is not Bessel-normalizable. This leaves us with an open question: Under what conditions can a Schauder basis be Bessel-normalizable? \qeddef
\end{remark}

We recall that a sequence $\set{x_n}$ is $\ell^2$-independent if $\sum c_nx_n=0$ implies that $c_n=0$ for all $n$ for any sequence $(c_n)\in \ell^2.$
It was shown by Casazza et al.\ in \cite[Thm 1.4]{CKST08} that every finitely linearly independent unit-norm Bessel sequence $\set{x_n}$ is a union of two sets each of which is $\ell^2$-independent.
We show next that every finitely linearly independent Bessel-normalizable sequence that is norm-bounded above also admits such a decomposition. 
\begin{proposition}
	Every finitely linearly independent Bessel-normalizable sequence that is norm-bounded above can be written as a union of two sets that are each $\ell^2$-independent.
\end{proposition}
\begin{proof}

	By \cite[Thm. 1.4]{CKST08}, $\bigset{\frac{x_{n}}{\norm{x_{n}}}}$ can be decomposed into $\bigset{\frac{x_{n}}{\norm{x_{n}}}}=\bigset{\frac{x_{n}^1}{\norm{x_{n}^1}}}\cup \bigset{\frac{x_{n}^2}{\norm{x_{n}^2}}}$ where $\bigset{\frac{x_{n}^i}{\norm{x_{n}^i}}}$ is a $\ell^2$-independent sequence for each $i$.
	Fix $1\leq i\leq 2$ and  let $(d_n)\in \ell^2 $ be such that $\sum\, d_n x_n^i=0$. Then we have $\sum \bigparen{d_n\,\norm{x^i_{n}}}\frac{x_n^i}{\norm{x^i_n}}=0$. Since $\bigparen{\,\norm{x_n}\,}$ is bounded above, $\bigparen{d_n\,\norm{x^i_{n}}}$ is in $\ell^2$. Hence, $d_n=0$ for all $n$, and consequently each sequence $\set{x_n^i}$ is a  $\ell^2$-independent sequence.
	\end{proof}

\subsection{Normalizable Bessel Sequences and Frames}
In this subsection, we study the normalizability of Bessel sequences and frames for $H$.

Bessel sequences are not necessarily frames. However, we show next that if $\set{x_n}$ is a frame, then $\nf$ is a Bessel sequence if and only if $\nf$ is a frame.
\begin{proposition}
	\label{equiv_nor_Be_fr}
	Assume that $\set{x_n}$ is a frame for $H$. Then $\set{x_n}$ is a normalizable frame if and only if $\set{x_n}$ is a normalizable Bessel sequence. 
\end{proposition}
\begin{proof}
	One direction is clear. For the other direction, assume that $\set{x_n}$ is a frame such that $\nf$ is a Bessel sequence. Let $S$ be the frame operator for $\set{x_n}$. Recall that  $\set{S^{-1/2}x_n}$ is a Parseval frame by Proposition $\ref{trans_frame_into_pars}$ and hence $\|S^{-1/2}x_n\|\leq 1$ for all $n$. Since $\|S^{-1/2}x_n\|\geq \|S^{1/2}\|^{-1}\|x_n\|$, it follows that \begin{equation*}
	\sum\, \Bigabs{\Bigip{x}{\frac{\|S^{1/2}\,\|x_n}{\|x_n\|}}}^2\Ge\sum \,\Bigabs{\Bigip{x}{\frac{x_n}{\|S^{-1/2}x_n\|}}}^2\Ge \sum\, |\ip{x}{x_n}|^2 \Ge A\,\|x\|^2,
	\end{equation*}
	where $A$ is a lower frame bound of $\set{x_n}$.
	Thus, $\sum \bigabs{\bigip{x}{\frac{x_n}{\|x_n\|}}}^2\Ge \frac{A}{\|S^{1/2}\|^2}\,\|x\|^2$. Consequently, $\nf$ is a frame. 
\end{proof}

Fix a sequence of scalars $c=(c_n)$ and define a diagonal operator $D_c:\text{Dom} (D_c)\rightarrow \ell^2$ by 
\begin{equation*}
D_c(d_n) = (c_nd_n), \qquad \text{for~} (d_n)\in \text{Dom}(D_c),
\end{equation*}
where $\text{Dom} (D_c)=\bigset{(d_n)~|~(c_nd_n)\in \ell^2}$. 
Let $\set{x_n}$ be a frame with associated analysis operator $C$. 
Then the analysis operator (not necessarily bounded) of $\set{c_nx_n}$ is $D_c\,\circ C$. 
Based on this perspective, Kutyniok et al.\ in \cite{KOPT13} studied frames which remain frames after being rescaled. The process of normalizing a frame is equivalent to rescaling the frame using the scalars $c_n=\frac{1}{\|x_n\|}$.
Combining Proposition \ref{equiv_nor_Be_fr} with \cite[Prop.\ 2.2]{KOPT13}, we obtain the following characterizations of normalizable frames.
We mention that the equivalence of (b) and (c) in Proposition \ref{diag_norm_fr} was given in  \cite[Prop.\ 2.2]{KOPT13}. For the sake of completeness, we include the proof here. For convenience, we say that a closed bounded linear operator $T$ between two Hilbert spaces that is injective and has closed range is Injective Closed Range or ICR (equivalently, $T$ is bounded and linear and there exists a constant $\delta>0$ such that $\|Tx\|\geq \delta \,\|x\|$ for all $x$).

\begin{proposition}
	\label{diag_norm_fr}
	Assume that $\set{x_n}$ is a frame for $H$ with associated analysis operator $C$. Let $\alpha=\bigparen{\frac{1}{\|x_n\|}}$. Then the following statements are equivalent.
	\begin{enumerate}
		\setlength\itemsep{0.6em}
		\item [\textup{(a)}] $\nf$ is a Bessel sequence in $H$.
		\item [\textup{(b)}] $\nf$ is a frame for $H$.
		\item [\textup{(c)}] Range($C$)$\,\subseteq \text{Dom}(D_\alpha)$ and $D_\alpha|_{\text{Range}(C)}$ is ICR.
		\item [\textup{(d)}] Range($C$)$\,\subseteq \text{Dom}(D_\alpha)$.
	\end{enumerate}
\end{proposition}

\begin{proof}
	(a) $\Rightarrow$ (b) By Proposition \ref{equiv_nor_Be_fr}.
	
	\medskip
	(b) $\Rightarrow$ (c) Assume that $\nf$ is a frame with associated analysis operator $\widetilde{C}$. Then for any $x\in H$ we have \begin{equation}
	\bigparen{\widetilde{C}x}_j\Eq\Bigip{x}{\frac{x_j}{\|x_j\|}}\Eq\bigparen{(D_\alpha \circ C)(x)}_j.
	\end{equation}
	Consequently, $\widetilde{C}=D_\alpha\circ C$. This shows that Range($C$) $\subseteq \text{Dom}(D_\alpha)$. Range($C$) is a closed subspace of $H$ because $\set{x_n}$ is a frame. It follows that $D_\alpha|_{\text{Range}(C)}$ is a closed, bounded linear operator between two Hilbert spaces. Let $A$ be a lower frame bound for $\nf$. Then for $x=Cy\in$ Range($C$)  we have \begin{equation}
	\|D_\alpha|_{\text{Range}(C)}(x)\|^2=\|(D_\alpha \circ C)(y)\|^2\Eq \norm{\widetilde{C}y}^2\Ge A\,\|y\|^2 \Ge A\, \|C\|^{-2}\,\|x\|^2.
	\end{equation}
	Thus, $D_\alpha|_{\text{Range}(C)}$ is ICR.
	\medskip
	
	(c) $\Rightarrow$ (d) This implication is clear.
	\medskip
	
	(d) $\Rightarrow$ (a) This implication follows by the Banach-Steinhaus Theorem. \qedhere
	
\end{proof}

The orthogonal projection of a frame onto a closed subspace is a frame for that closed subspace.
Combining this fact with Theorem \ref{fr_equi_orth_basis_el2}, Lemma \ref{topo_preserve_normal}, and Proposition \ref{equiv_nor_Be_fr}, we obtain another characterization of normalizable frames.
\begin{proposition}
	\label{norm_orp_basis}
	Assume that $\set{x_n}$ is a frame for $H$ and let $S$ be the range of the analysis operator of $\set{x_n}$. Then $\set{x_n}$ is a normalizable frame for $H$ if and only if $\set{P_{S}\delta_n}$ is a normalizable Bessel sequence in $S$, where $\set{\delta_n}$ is the standard basis for $\ell^2$. \qeddef
\end{proposition}

\begin{corollary}
	Let $\set{x_n}$ and $\set{y_n}$ be two frames for $H$ with associated analysis operators $C$ and $\widetilde{C}$, respectively. Assume that $\text{Range}(C)\Eq \text{Range}(\widetilde{C})$. Then $\set{x_n}$ is a normalizable frame if and only if $\set{y_n}$ is a normalizable frame. \qeddef
\end{corollary}

We saw in Example \ref{nbb_normalizable} that every frame that is norm-bounded below is a normalizable frame, but the converse does not hold in general (see Example \ref{non_nor_frames}). The next result shows that every normalizable frame must contain a subsequence that is a frame that is norm-bounded below. 


\begin{theorem}
	\label{norma_cont_nbb_frame}
	If $\set{x_n}$ is a normalizable frame for $H$, then it falls into exactly one of the following three categories.
	\begin{enumerate}
		\setlength\itemsep{0.25em}
		\item [\textup{(a)}] $\set{x_n}$ is a frame that is norm-bounded below.		
		\item [\textup{(b)}] There exists a $\delta>0$ such that $\set{x_n~|~\|x_n\|\geq \delta }$ is a frame and $\set{x_n~|~\|x_n\|< \delta }$ is not empty but does not satisfy a lower frame condition.		
		\item [\textup{(c)}] There exists a sequence of positive scalars $(\delta_k)$ that decreases to $0$ such that $\mathcal{F}_k=\set{x_n~|~\delta_{k+1}\leq\|x_n\|< \delta_k}$ is a frame for every $k$. In particular, if $A_k$ and $B_k$ are optimal lower and upper frame bounds for $\mathcal{F}_k$, respectively, then $(A_k)\in \ell^1$ and $\lim\limits_{k\rightarrow \infty} B_k=0.$
	\end{enumerate}
\end{theorem}
\begin{proof}
	We first show that $\set{x_n}$ must contain a frame that is norm-bounded below.
	Let $C$ be an upper frame bound for $\nf$ and let $A$ be a lower frame bound for $\set{x_n}$. Suppose to the contrary that $\set{x_n}$ does not contain any frame that is norm-bounded below.  Fix $0<\delta< \sup \|x_n\|$. Note that $\set{x_n~|~\|x_n\|\geq \delta }$ is a Bessel sequence. 
	Since $\set{x_n~|~\|x_n\|\geq \delta }$ is not a frame, there exists some nonzero element $x_\delta \in H$ such that
	\begin{equation}
	\label{lower_frame_vio}
	\sum_{\|x_n\|\Ge \delta}|\ip{x_\delta}{x_n}|^2\,<\,\frac{A}{2}\,\|x_\delta\|^2.
	\end{equation}
	Equation (\ref{lower_frame_vio}) implies that 
	\begin{equation*}
	\sum_{\|x_n\|\,<\,\delta}|\ip{x_\delta}{x_n}|^2\Eq \sumli |\ip{x_\delta}{x_n}|^2\,-\sum_{\|x_n\|\Ge \delta}|\ip{x_\delta}{x_n}|^2\,\Ge\, \frac{A}{2}\,\|x_\delta\|^2.
	\end{equation*}
	Therefore,
	
	\begin{equation*}
	C\,\|x_\delta\|^2\,\Ge\, \sum_{\|x_n\|\,<\,\delta}\,\Bigabs{\Bigip{x}{\frac{x_n}{\norm{x_n}}}}^2\Ge  \frac{1}{\delta^2}\sum_{\|x_n\|\,<\,\delta}|\ip{x}{x_n}|^2 \Ge \frac{A}{2\delta^2}\,\|x_{\delta}\|^2.
	\end{equation*}
	Since $x_\delta$ is nonzero, we conclude that $C>\frac{A}{2\delta^2}$. But $\delta$ can be arbitrarily small, so we obtain a contradiction.
	
	Now, let $\set{x_{n_k}}$ be a frame that is norm-bounded below. 
	Choose  $ \delta_1>\sup \|x_n\|$ and set $\delta_2=\inf \|x_{n_k}\|>0$. Then $\set{x_n\,|\,\delta_2\leq \|x_n\|< \delta_1}$ is a frame since every subsequence of a Bessel sequence is a Bessel sequence and 
	$\sum_{\,\delta_2\leq \|x_n\|< \delta_1\,}|\ip{x}{x_n}|^2\geq \sum |\ip{x}{x_{n_k}}|^2$. If $\set{x_n\,|\,\|x_n\|<\delta_2}$ is empty, then $\set{x_n}$ belongs to category (a).  If $\set{x_n\,|\,\|x_n\|<\delta_2}$ is non-empty and satisfies a lower frame condition, then it is a normalizable frame by Proposition \ref{equiv_nor_Be_fr}. Arguing similarly to the process of finding $\delta_2$, then we can find $\delta_3$ such that $\set{x_n\,|\,\delta_3\leq\|x_n\|<\delta_2}$ is a frame. If $\set{x_n\,|\,\|x_n\|<\delta_3}$ is empty, then $\set{x_n}$ belongs to category (a) since $\set{x_n\,|\, \|x_n\|\geq \delta_3}$ is the union of two frames, and hence is a frame. Assume that 
	$\set{x_n\,|\, \|x_n\|< \delta_3}$ does not satisfy the lower frame condition, then $\set{x_n}$ falls into category (b). Otherwise, iterating this process, we either stop at some stage (ending in category (a) or (b)) or the iteration never ends (putting us in category (c)). 
	
	It remains to show that $(A_k)\in \ell^1$ and $\lim\limits_{k\rightarrow \infty} B_k=0$ if $\set{x_n}$ falls into category (c). 
	Since $$\Bigparen{\sum A_k}\norm{x}^2 \Le \sum_k \sum_{x_n\in \mathcal{F}_k} |\ip{x}{x_n}|^2\Le \sum\, |\ip{x}{x_n}|^2,$$
	we have that $(A_k)\in \ell^1$. Suppose to the contrary that the $B_k$ do not converge to 0 as $k\rightarrow \infty$. Then there exists a constant $\alpha >0$ such that for any $k\in \N$ there exists $n_k\geq k$ for which $B_{n_k}\geq \alpha.$ For each $n_k$, let $y_{n_k}\in H$ with $\norm{y_{n_k}}=1$ be such that $\sum_{x_n\in \mathcal{F}_{n_k}} |\ip{y_{n_k}}{x_n}|^2\geq \frac{B_{n_k}}{2}.$ However, this implies that 
	$$\sum_{x_n\in \mathcal{F}_{n_k}} |\ip{y_{n_k}}{\frac{x_n}{\|x_n\|}}|^2\Ge \frac{\alpha}{2\delta_{n_k}}\rightarrow \infty$$ as $k\rightarrow \infty$, which contradicts that $\nf$ is a Bessel sequence.
\end{proof}

Based on our observations, we propose the following conjecture.
\begin{conjecture}
	There does not exist a normalizable frame belonging to category (c) in Theorem \ref{norma_cont_nbb_frame}. \qeddef
\end{conjecture}

\subsection{Sequences That Are Not Frame-Normalizable}
Theorem \ref{norma_cont_nbb_frame} shows that a normalizable frame must contain a frame that is norm-bounded below. Consequently, we have the following sufficient condition for a frame to not be a normalizable frame.

\begin{corollary}
	Every frame that does not contain a frame that is norm-bounded below is not a normalizable frame. \qeddef
\end{corollary}

We provide some examples in the following.

One way to construct a frame that does not contain a frame that is norm-bounded below is to duplicate elements infinitely many times but rescale their coefficients so that we can obtain a frame. 
\begin{example}
	\label{non_nor_frames}
	Let $\set{e_n}$ be an orthonormal basis for $H$. Then
	\begin{equation*} \set{x_n}=\Bigset{e_1,\,\frac{1}{\sqrt{2}}e_2,\,\frac{1}{\sqrt{2}}e_2,\,\frac{1}{\sqrt{3}}e_3,\,\frac{1}{\sqrt{3}}e_3,\,\frac{1}{\sqrt{3}}e_3,\,\dots} \end{equation*} is a Parseval frame for $H$. However, $\nf$  is not a Bessel sequence.\qeddef
\end{example}

An exact Bessel sequence need not be normalizable.

\begin{example}
	\label{non_nor_ex_Be_se}
	Let $\set{e_n}$ be an orthonormal basis for $H$. Define $x_n=\frac{1}{n}e_1+\frac{1}{n}e_{n+1}$ for $n\in\mathbb{N}$. We have that for any $x\in H$,
	\begin{equation*}
	\sum \,	|\ip{x}{x_n}|^2\Eq \sum\frac{1}{n^2}\,|\ip{x}{e_1}+\ip{x}{e_{n+1}}|^2\Le \Bigparen{\,\sum \frac{2}{n^2}\,}\|x\|^2.
	\end{equation*}
	Consequently, $\set{x_n}$ is a Bessel sequence.
	
	It is clear that the sequence $a_n=ne_{n+1}$ is biorthogonal to $\set{x_n}$, so $\set{x_n}$ is minimal. If $x=\sum c_ke_k\in H$ is such that  $\ip{x}{x_n}=0$ for all $n$, then $c_1=-c_k$ for all $k\geq 2$. Therefore, $x=0$, so $\set{x_n}$ is exact.
	However, $\frac{x_n}{\|x_n\|}=\frac{\sqrt{2}}{2}e_1+\frac{\sqrt{2}}{2}e_{n+1}$ for each $n$. Since $\sum\, \bigabs{\bigip{e_1}{\frac{x_n}{\|x_n\|}}}^2=\infty$, it follows that $\set{x_n}$ is not a normalizable Bessel sequence.\qeddef
\end{example}

A frame of the form $\set{A^n x}_{n\geq 0}$, where $x\in H$ is fixed, cannot be a normalizable frame if $A\colon H\rightarrow H$ is a self-adjoint operator \cite[Thm 6.2]{ACMCP16}. In the same paper, conditions for such sequences to be frames when $A$ is a self-adjoint operator were also provided. We summarize their results as follows (see also \cite[Example 4]{ACMCP16}).   

\begin{theorem} [Non-normalizable frames generated by iterations of self-adjoint operators]
	\label{non_norm_frame_self_ad_oper}
	Let $A\colon H\rightarrow H$ be a self-adjoint operator on an infinite-dimensional, separable Hilbert space $H$, and fix $x\in H$. Assume that:
	\begin{enumerate}
		\setlength\itemsep{0.3em}
		\item[\textup{(a)}]$A=\sum_j \lambda_jP_j$, where the $P_j$'s are rank-one orthogonal projections,
		\item[\textup{(b)}]$|\lambda_k|<1$ for all $k$,
		\item[\textup{(c)}]$\lim\limits_{k\rightarrow \infty} |\lambda_k|=1$,
		\item[\textup{(d)}]$\displaystyle \inf_n \prod_{k\neq n} \frac{|\lambda_n-\lambda_k|}{|1-\overline{\lambda_n}\lambda_k|} \geq \delta$ (\emph{Carleson's condition}) for some $\delta>0$, and
		\item[\textup{(e)}]  there exist constants $C_1,C_2$ such that $0<C_1\leq \frac{\|P_j x\|}{\sqrt{1-|\lambda_k|^2}}\leq C_2<\infty$ for every $k,j$.
	\end{enumerate} 
	Then $\set{A^n x}_{n\geq 0}$ is a frame but it is not normalizable. \qeddef
\end{theorem}

\begin{remark}
	We remark that $\set{x_n+y_n}$ needs not be a Bessel-normalizable sequence even if $\set{x_n}$ and $\set{y_n}$ are Bessel-normalizable sequences. This was first proved for frames by Casazza (\cite{Cas98}) and later extended to Bessel sequences that every Bessel sequence in a complex Hilbert space is a linear combination of two Riesz bases (see \cite{DM15} or \cite{HY23}). By Example \ref{non_nor_frames}, let $\set{z_n}$ be a Bessel sequence that is not Bessel-normalizable. Then there exist two Riesz bases $\set{x_n}$ and $\set{y_n}$ such that $z_n=x_n+y_n$ for all $n$. Our remark then follows from that every Riesz basis is a frame that is norm-bounded below and hence is frame-normalizable. \qeddef
\end{remark}

 Example \ref{non_nor_frames} constructed a frame $\set{x_n}$ for which $\|x_n\|$ converges to 0 as $n \rightarrow \infty$ that is not Bessel-normalizable . We generalize this example by proving that every frame for which $\|x_n\|$ converges to 0 as $n\rightarrow \infty$ is not Bessel-normalizable.
 
The following lemma follows from straightforward computations, which we omit here.
\begin{lemma}
	\label{infi_zero_lemma}
	Let $\set{x_n}$ be a sequence in $H$ and let $(c_n)$ be a sequence of scalars for which $c_n\neq 0$ for all $n$. Then the following statements hold.
	\begin{enumerate}
		\item[\textup{(a)}] If $c_n$ converges to $0$ and $\bigset{\frac{x_n}{c_n}}$ is a Bessel sequence, then $\set{x_n}$ is a Bessel sequence.
		\item[\textup{(b)}] If $c_n$ converges to $\infty$, then $\set{x_n}$ satisfies a lower frame condition if $\bigset{\frac{x_n}{c_n}}$ satisfies a lower frame condition. Moreover, $\bigset{\frac{x_n}{c_n}}$ is a Bessel sequence if $\set{x_n}$ is a Bessel sequence. \qeddef
	\end{enumerate} 
\end{lemma}

\begin{theorem}
	\label{non_norma_decrea_norm}
	Let $\set{x_n}\subseteq H$ be a sequence for which $\norm{x_n}$ converges to $0$ as $n\rightarrow \infty$. Assume that there exists an infinite-dimensional closed subspace $M$ of $H$ such that the sequence $\set{P_Mx_{n}}$ satisfies a lower frame condition for $M$. Then $\set{x_n}$ is not Bessel-normalizable.
\end{theorem} 
\begin{proof}
	Suppose to the contrary that $\set{x_n}$ is Bessel-normalizable. By Lemma \ref{infi_zero_lemma}, $\set{x_n}$ is also a Bessel sequence and hence $\set{P_Mx_{n}}$ is a Bessel sequence in $M$. 
	Since $\set{P_Mx_{n}}$ satisfies a lower frame condition for $M$, it is a frame for $M$. Let $T\colon M\rightarrow M$ be the frame operator for $\set{P_Mx_{n}}$. Then $\set{T^{-1/2}P_Mx_{n}}$ is a Parseval frame for $M$ by Proposition \ref{trans_frame_into_pars}. By the Bessel-normalizability of $\set{x_n}$, we see that $\Bigset{\frac{P_Mx_{n}}{\|x_{n}\|}}$ is a Bessel sequence in $M$.
	Consequently, $\Bigset{T^{-1/2}\frac{P_Mx_{n}}{\|x_{n}\|}}$ is a Bessel sequence in $M$. Let $C$ be an upper frame bound for $\Bigset{T^{-1/2}\frac{P_Mx_{n}}{\|x_{n}\|}}$ and choose $K$ and $n_0$ large enough that $\frac{3K}{4}>C$ and $\frac{1}{\|x_{n}\|^2}\geq K$ for $n\geq n_0$. Let $J=\set{n}_{n=n_0}^\infty$. Then by Proposition \ref{iden_parseval}, we have
	\begin{align}
	\label{first_eq}
	\begin{split}
	C\|x\|^2 &\Ge \sum_{n=n_0}^\infty\, \Bigabs{\Bigip{x}{\frac{T^{-1/2}P_Mx_{n}}{\|x_{n}\|}}}^2\\ 
	&\Ge K \sum_{n=n_0}^\infty \,\Bigabs{\Bigip{x}{T^{-1/2}P_Mx_{n}}}^2\\
	&\Ge K\,\Bigparen{\,\frac{3}{4}\|x\|^2\,-\,\Bignorm{\sum_{n\in J^c}\,\Bigip{x}{T^{-1/2}P_Mx_{n}}\,T^{-1/2}P_Mx_{n}}^2\,}.	
	\end{split}		
	\end{align}
	Since $S=\clspan{\set{T^{-1/2}P_Mx_{n}}_{n\in J^c}}$ is finite-dimensional, there exists some nonzero element $x_0 \in M\cap  S^\perp$. Inserting $x_0$ into inequality (\ref{first_eq}), we see that $C\|x_0\|^2\geq \frac{3K}{4}\|x_0\|^2 > C\|x_0\|^2,$ which is a contradiction.
\end{proof}

\begin{corollary}
	\label{non_norma_frame_with_norm_conv_zero}
	Let $\set{x_n}$ be a frame for $H$. If $\|x_n\|$ converges to 0 as $n\rightarrow \infty$, then $\set{x_n}$ is not Bessel-normalizable. \qeddef
\end{corollary} 

By a similar argument, we obtain a lower-normalizable version of Theorem \ref{non_norma_decrea_norm} as follows.
\begin{theorem}
	\label{non_norma_increa_norm}
		Let $\set{x_n}\subseteq H$ be a sequence for which $\norm{x_n}$ converges to $\infty$ as $n \rightarrow \infty$. Assume that there exists an infinite-dimensional closed subspace $M$ of $H$ such that the sequence $\set{P_Mx_{n}}$ is a Bessel sequence in $M$. Then $\set{x_n}$ is not lower-normalizable.
\end{theorem} 
\begin{proof}
	Suppose that $\set{x_n}$ is lower-normalizable. Since $\set{P_Mx_{n}}$ is a Bessel sequence, $\Bigset{\frac{P_Mx_{n}}{\|x_n\|}}$ is also a Bessel sequence in $M$ by Lemma \ref{infi_zero_lemma}.
	By the lower-normalizability of $\set{x_n}$, we see that $\Bigset{\frac{P_Mx_{n}}{\|x_n\|}}$ is a frame for $M$. 
	Let $T\colon M\rightarrow M$ be the frame operator for $\Bigset{\frac{P_Mx_{n}}{\|x_n\|}}$. Then we have that  $\Bigset{T^{-1/2}\frac{P_Mx_{n}}{\|x_n\|}}$ is a Parseval frame for $M$. Let $C$ be an upper frame bound for $\set{T^{-1/2}P_Mx_{n}}$ and choose $K$ and $n_0$ large enough that $\frac{3}{4}>\frac{C}{K}$ and $\|x_{n}\|^2\geq K$ for $n\geq n_0$. Let $J=\set{n}_{n=n_0}^\infty$. We compute that 
	\begin{align}
	\label{second_eq}
	\begin{split}
	\frac{C}{K}\|x\|^2 &\Ge\, \frac{1}{K}\sum_{n=n_0}^\infty\, \Bigabs{\Bigip{x}{T^{-1/2}P_Mx_{n}}}^2\\ 
	&\Ge\,  \sum_{n=n_0}^\infty \,\Bigabs{\Bigip{x}{T^{-1/2} \frac{P_Mx_{n}}{\|x_n\|} }}^2\\
	&\Ge\Bigparen{\,\frac{3}{4}\|x\|^2\,-\,\Bignorm{\sum_{n\in J^c}\,\Bigip{x}{T^{-1/2}\frac{P_Mx_{n}}{\|x_n\|}}\,T^{-1/2}\frac{P_Mx_{n}}{\|x_n\|}}^2\,}.	
	\end{split}		
	\end{align}
	where the last inequality follows by Proposition \ref{iden_parseval}.
	Since $S=\clspan{\Bigset{T^{-1/2}\frac{P_Mx_{n}}{\|x_n\|}}_{n\in J^c}}$ is finite-dimensional, arguing similarly to Theorem \ref{non_norma_decrea_norm}, we see that $\frac{C}{K}\Ge \frac{3}{4}$, which is a contradiction.
\end{proof}

\section{Perturbation of frame-normalizable sequences}
\label{pert_norm_frames}
In this section, we present some ways that allow us to perturb frame-normalizable sequences without breaking the frame-normalizablility.

Topological isomorphisms send frames to frames. Even more, topological isomorphisms preserve frame-normalizability and Bessel-normalizability.
\begin{lemma}
	\label{topo_preserve_normal}
	Let $H$ and $K$ be separable Hilbert spaces and let $T:H\rightarrow K$ be a topological isomorphism. Then $\nf$ is a frame (\textit{resp. Bessel sequence}) for $H$ if and only if $\bigset{\frac{T(x_n)}{\|T(x_n)\|}}$ is a frame for K (\textit{resp. Bessel sequence}).
\end{lemma}
\begin{proof}
	The proof follows from the standard inequalities $\norm{T^{-1}}^{-1}\norm{x}\leq \norm{Tx} \leq \norm{T}\norm{x}$ and a similar inequality for $T^*$.
\end{proof}

The classic Paley-Wiener perturbation theorem states that if $\set{x_n}$ is a Schauder basis for a Banach space $X$, $\set{y_n}\subseteq X$, and there exists $0\leq\lambda<1$  such that $$\Bignorm{\sum_{n=1}^N c_n(x_n-y_n)}\Le \lambda\, \Bignorm{\sum_{n=1}^N c_nx_n}$$
for any $N\in \mathbb{N}$ and scalars $c_1,...,c_N$, then $\set{y_n}$ is a Schauder basis for $X$. This result was later extended to frames in Hilbert spaces by Casazza and Christensen in \cite{CC97}. 
\begin{theorem}[\cite{CC97}]
	\label{pertu_chris_frames}
	Let $\set{x_n}$ be a frame for $H$ with frame bounds $A,B$. Let $\set{y_n}\subseteq H$ and assume there exist constants $\lambda, \mu,\nu \geq 0$ such that $\max\bigset{\lambda+\frac{\mu}{\sqrt{A}},\nu}<1$. Assume that for all $N\in \mathbb{N}$ and scalars $c_1,...,c_N$,
	\begin{equation}
	\label{perturb_ineq}
	\Bignorm{\sum_{n=1}^N c_n(x_n-y_n)}\Le \lambda \,\Bignorm{\sum_{n=1}^N c_nx_n}+\mu\, \biggparen{\,\sum_{n=1}^N |c_n|^2\,}^{\!1/2}+\,\nu\,\Bignorm{\sum_{n=1}^N c_ny_n}.
	\end{equation}
 Then $\set{y_n}$ is a frame for $H$ with frame bounds $A\,\Bigparen{1-\Bigparen{\frac{\lambda+\nu+\frac{\mu}{\sqrt{A}}}{1+\nu}}}^2,\, B\,\Bigparen{1+\Bigparen{\frac{\lambda+\nu+\frac{\mu}{\sqrt{B}}}{1-\nu}}}^2$.\qeddef
\end{theorem}

In addition to Theorem \ref{pertu_chris_frames}, we refer to \cite{CC98} and \cite{CH97} for other type of perturbation theorems for frames. 

We first observe that the criterion $\|\sum_{n=1}^N c_n(x_n-y_n)\|\Le\mu\, (\sum_{n=1}^N |c_n|^2)^{1/2}$ in Theorem \ref{pertu_chris_frames} may not preserve frame-normalizability. By the Cauchy–Bunyakovsky–Schwarz inequality, we know that $\|\sum_{n=1}^N c_n(x_n-y_n)\|\Le\mu\, (\sum_{n=1}^N |c_n|^2)^{1/2}$ if $\bignorm{ \,\bigparen{\,\|x_n-y_n\|\,}\,}_{\ell^2}<\mu$. Choosing the coefficients properly, we can construct two sequences $\set{x_n}$ and $\set{y_n}$ for which the $\ell^2$-norm of $(\,\|x_n-y_n\|\,)$ is small enough, but $\bigset{\frac{y_n}{\|y_n\|}}$ contains infinitely many copies of some element. 

The following lemma follows by straighforward computation, which we omit here.
\begin{lemma}
	\label{norm_ratio_equi_normal_fr}
	Let $\set{x_n}\subseteq H$ be a sequence and let $(c_n)$ be a sequence of scalars. Assume that there exist constants $0< M\leq L$ such that $M\leq  \frac{\|x_n\|}{|c_n|}\leq L$ for every $n$. Then $\nf$ is a frame if and only if $\bigset{\frac{x_n}{c_n}}$ is a frame.\qeddef
\end{lemma}

Our perturbation theorem for frame-normalizable sequences is as follows.
\begin{proposition}
	\label{pertu_thm_norma_fr}
	Let $\set{x_n}\subseteq H$ be a frame-normalizable sequence and let $A$ be a lower frame bound for $\nf$.
	Assume that $\set{y_n}\subseteq H$ is a sequence not containing any zero elements that satisfies one of the following conditions for every $N\in \mathbb{N}$ and scalars $c_1,...,c_N$,
	\begin{enumerate}
		\setlength\itemsep{0.3em}
		\item [\textup{(a)}] There exist $\lambda, \nu \in [0,1)$ such that 
		\begin{equation}
		\label{perturb_ineq2}
		\Bignorm{\sum_{n=1}^N c_n(x_n-y_n)}\Le \lambda \,\Bignorm{\sum_{n=1}^N c_nx_n} \,+ \,\nu\,\Bignorm{\sum_{n=1}^N c_ny_n}.
		\end{equation}
		
		\item [\textup{(b)}] There exists some $0\leq K<\min\set{1,\sqrt{A}\,}$ such that \begin{equation} \label{perturb_ineq3}\Bignorm{\sum_{n=1}^N \frac{c_n}{\norm{x_n}}(x_n-y_n)} \,\Le K \,\biggparen{\,\sum_{n=1}^N |c_n|^2\,}^{\!1/2}. \end{equation}
	
		\item [\textup{(c)}] There exists some $0\leq K<\frac{\sqrt{A}}{1+\sqrt{A}}$ such that \begin{equation} \label{perturb_ineq4} \Bignorm{\sum_{n=1}^N \frac{c_n}{\norm{y_n}}(x_n-y_n)} \,\Le K \,\biggparen{\,\sum_{n=1}^N |c_n|^2\,}^{\!1/2}. \end{equation}
		
	\end{enumerate}
		Then $\set{y_n}$ is a frame-normalizable sequence. 
\end{proposition}

\begin{proof}
	By hypothesis, we have that $\nf$ is a frame for $H$.
	\medskip
	
	(a) By considering $c_n=\frac{c_n}{\norm{x_n}}$ in inequality (\ref{perturb_ineq2}), we see that for any $N\in \N$ and scalars $c_1,...,c_N$, $$\Bignorm{\sum_{n=1}^N \frac{c_n}{\|x_n\|}(x_n-y_n)}\Le \lambda \,\Bignorm{\sum_{n=1}^N c_n\frac{x_n}{\|x_n\|}} \,+\, \nu\,\Bignorm{\sum_{n=1}^N c_n\frac{y_n}{\|x_n\|}}.$$
	By Theorem \ref{pertu_chris_frames}, $\bigset{\frac{y_n}{\|x_n\|}}$ is a frame for $H$. Fix $1\leq n\leq N$ and set $c_n=1$ and $c_i=0$ for $i\neq n$ in inequality (\ref{perturb_ineq2}), we have 
	$\norm{ x_n-y_n}\Le \lambda \,\norm{x_n} \,+ \,\nu\,\norm{y_n}. $
	By the Triangle Inequality, we obtain \begin{equation} 
	(1-\nu)\,\|y_n\|\leq (1+\lambda)\,\|x_n\| \quad \text{and}  \quad (1-\lambda)\,\|x_n\|\leq (1+\nu)\,\|y_n\|.
	\end{equation}
	Therefore, $\frac{1-\lambda}{1+\nu} \Le \frac{\|y_n\|}{\|x_n\|}\Le \frac{1+\lambda}{1-\nu} $ for all $n$. Hence,  $\bigset{\frac{y_n}{\|y_n\|}}$ is a frame for $H$ by Lemma \ref{norm_ratio_equi_normal_fr}.
	
	\medskip
	(b) 
	By Theorem \ref{pertu_chris_frames}, we see that $\bigset{\frac{y_n}{\norm{x_n}}}$ is a frame for $H$.  Fix $1\leq n\leq N$ and set $c_n=1$ and $c_i=0$ for $i\neq n$ in inequality (\ref{perturb_ineq3}), we have $\norm{x_n-y_n}\Le K\norm{x_n}$ for any $1\leq n\leq N$. Using the Triangle Inequality, we obtain $$(1-K)\,\norm{x_n}\Le \norm{y_n} \Le (1+K)\,\norm{x_n}.$$
    It follows that $\bigset{\frac{y_n}{\norm{y_n}}}$ is a frame for $H$ by Lemma \ref{norm_ratio_equi_normal_fr}.
	
	\medskip
	(c) Similar to part (b), fix $1\leq n\leq N$ and set $c_n=1$ and $c_i=0$ for $i\neq n$ in inequality (\ref{perturb_ineq4}), we see that $(1-K)\,\norm{y_n}\Le \norm{x_n} \Le (1+K)\,\norm{y_n}$. Equivalently, $(1-K)\Le \frac{\norm{x_n}}{\norm{y_n}} \Le (1+K)$. Since $\nf$ is a frame for $H$ with a lower frame bound $A$, it follows that $\bigset{\frac{x_n}{\norm{y_n}}}$ is a frame with a lower frame bound $(1-K)^2A.$ Because $0<K<\frac{\sqrt{A}}{1+\sqrt{A}}$, we obtain $K(1+\sqrt{A})<\sqrt{A}$ and hence $\frac{K}{(1-K)\sqrt{A}}<1$. By Theorem \ref{pertu_chris_frames}, $\bigset{\frac{y_n}{\norm{y_n}}}$ is a frame for $H$.
\end{proof}
	

\begin{remark} 
		(a)
		We mention that Proposition \ref{pertu_thm_norma_fr} also holds for Bessel-normalizable sequences by using the Triangle Inequality and \cite[Thm. 7.4]{Hei11}.\medskip
		
        (b)
		The conditions $\lambda\in[0,1)$ in statement (a) and $0\leq K<\min\set{1,\sqrt{A}\,} $ in statement (b) are tight. To show this, we use the idea of \cite[Example 5.4.6]{Chr16}. Let $\set{e_n}$ be an orthonormal basis. Define $x_n=e_n$ and $y_n=e_n+e_{n+1}$ for all $n$. Then $\nf$ is an orthonormal basis with a lower frame bound $1$. In particular, for all $N\in \mathbb{N}$ and scalars $c_1,...c_N,$ we have 
		$$\Bignorm{\sum_{n=1}^N c_n(x_n-y_n)}\Eq \Bignorm{\sum_{n=1}^N \frac{c_n}{\norm{x_n}}(x_n-y_n)}\Eq \biggparen{\,\sum_{n=1}^N |c_n|^2\,}^{\!1/2}\Eq \Bignorm{\sum_{n=1}^N c_nx_n}.$$
	 However, since $e_1$ cannot be expressed in the form of $\sum c_n\frac{y_n}{\norm{y_n}}$ for any choice of coefficients of $c_n$, it follows that $\bigset{\frac{y_n}{\norm{y_n}}}$ is not a frame for $H$.\\

		(c) For any $\mu>0$, the condition $\|\sum_{n=1}^N c_n(x_n-y_n)\|\Le \mu \,(\,\sum_{n=1}^N |c_n|^2\,)^{1/2}$ does not necessarily preserve frame-normalizability even if the two sequences involved in the perturbation are frames for $H$. Let $\set{e_n}$ be an orthonormal basis for $H$ and let $(d_n)$ be a sequence of nonzero scalars which will be chosen later.
		Define $\set{x_n}$ to be $x_{2n-1}=e_n$ and $x_{2n}=d_n e_n$ for all $n$. On the other hand, set $\set{y_n}$ to be $y_{2n-1}=e_n$ and $y_{2n}=d_n e_1$ for all $n$. We then see that $\set{x_n}$ and $\set{y_n}$ are frames if $(d_n)\in \ell^2$. Moreover, $$\Bignorm{\sum_{n=1}^N c_n(x_n-y_n)}\Le \biggparen{\,\sum_{n=1}^\infty 2d_n^{\,2}\,}^{\!1/2}\biggparen{\,\sum_{n=1}^N |c_n|^2\,}^{\!1/2}.$$
		We can choose $(d_n)$ to be small enough so that $(\,\sum_{n=1}^\infty 2d_n^{\,2}\,)^{1/2}$ is smaller than $\mu$. However, $\set{y_n}$ is not a normalizable frame since $\bigset{\frac{y_n}{\|y_n\|}}$ contains infinitely many copies of $e_1$.\\

		(d) The counterexample given in (c) also serves as a counterexample that compactness of the linear operator $T\colon \ell^2 \rightarrow H$ defined by $T(c_n)=\sum_{n=1}^\infty c_n(x_n-y_n)$ does not preserve frame-normalizability.\\

		(e) Orthogonal projection does not necessarily preserve Bessel-normalizability. By Example \ref{non_nor_frames}, there exists a frame $\set{x_n}$ that is not Bessel-normalizable. By Theorem \ref{fr_equi_orth_basis_el2}, every frame is topologically isomorphic to the orthogonal projection of the standard basis of $\ell^2$ onto some closed subspace $S$ of $\ell^2$. Orthonormal bases are normalizable frames. By Lemma \ref{topo_preserve_normal}, we obtain an example that shows that an orthogonal projection of an orthonormal basis onto a closed subspace of $H$ needs not be a Bessel-normalizable sequence in the closed subspace.		\qeddef
\end{remark}

\section{Applications}
\label{BS_conjecture}
\subsection{The Balazs-Stoeva Conjecture}

We show that the Balazs-Stoeva conjecture holds for Bessel-normalizable sequences in this section.
The main tool is Orlicz's Theorem, which was proved in \cite{Or33} (for a textbook description see \cite[Thm. 3.16]{Hei11}). 
\begin{theorem}[Orlicz's Theorem]
	\label{Orlicz_thm}
	If $\set{x_n}$ is a sequence in $H$ for which $\sum x_n$ converges unconditionally, then $\sum \|x_n\|^2<\infty.$ \qeddef
\end{theorem}

\begin{definition}
Let $\mathcal{X}=\set{x_n} \subseteq H$ and $\mathcal{Y}=\set{y_n} \subseteq H$ be two sequences, and let $\mathcal{M}=(m_n)$ be a sequence of scalars. Then the \emph{multiplier} associated with $\mathcal{X}, \mathcal{Y},$ and $\mathcal{M}$ is the operator on $H$ defined by $\mathcal{M}_{\mathcal{X},\,\mathcal{Y}}\,(x)=\sum m_n\,\ip{x}{y_n}\,x_n.$  \qeddef
\end{definition}

The Balazs-Stoeva conjecture states that $\mathcal{M}_{\mathcal{X},\,\mathcal{Y}}\,(x)$ converges unconditionally for all $x\in H$ if and only if 
there exist sequences of scalars $(c_n)$ and $(d_n)$ with $c_n\overline{d_n}=m_n$ such that $\set{c_nx_n}$ and $\set{d_ny_n}$ are Bessel sequences. 

The direction from left to right of the following result was proved in \cite[Propostion 3.2]{BS11}. For the sake of completeness, we include the short proof here.
\begin{proposition}
	\label{Orlicz_coro_BS_conjecture}
	Let $\mathcal{X}=\set{x_n}$ be a Bessel-normalizable sequence in $H$. Assume that $\mathcal{Y}=\set{y_n}$ is a sequence in $H$ and $\mathcal{M}=(m_n)$ is a sequence of scalars. Then $\sum m_n \,\ip{x}{y_n}\,x_n$ converges unconditionally for all $x\in H$ if and only if $\set{m_n\,\|x_n\|\,y_n}$ is a Bessel sequence.
	
	As a consequence, since $m_n=m_n\,\|x_n\| \cdot \frac{1}{\|x_n\|}$, the Balazs-Stoeva conjecture holds for the multiplier  $\mathcal{M}_{\mathcal{X},\mathcal{Y}}$. 
	
\end{proposition}
\begin{proof}
	($\Rightarrow$)
	By Orlicz's Theorem and the Banach-Steinhaus Theorem, we see that $\set{\overline{m_n}\,\|x_n\|\,y_n}$ is a Bessel sequence.
	
	\medskip
	($\Leftarrow$) This direction follows from our hypothesis and \cite[Thm.\ 7.2]{Hei11}.
\end{proof}

\begin{corollary}
	\label{BS_conjecture_p_norma}
	Let $\mathcal{X}=\set{x_n}$ be a sequence in $H$ that is norm-bounded above. Assume that $\mathcal{Y}=\set{y_n}\subseteq H$ and $\mathcal{M}=(m_n)$ is a sequence of scalars. If $\bigset{\frac{x_n}{\|x_n\|^p}}$ is a Bessel sequence for some $1\leq p < \infty$, then the Balazs-Stoeva conjecture holds for the multiplier  $\mathcal{M}_{\mathcal{X},\,\mathcal{Y}}$. Specifically, $m_n=m_n\,\|x_n\|^p \cdot \frac{1}{\|x_n\|^p}$ and $\set{m_n\,\|x_n\|^p\, y_n}$ and $\bigset{\frac{x_n}{\|x_n\|^p}}$ are Bessel sequences.
	
\end{corollary}
\begin{proof}
	It suffices to show that $\set{\overline{m_n}\,\|x_n\|^p \,y_n}$ is a Bessel sequence if $\sum m_n \,\ip{x}{y_n}\,x_n$ converges unconditionally for all $x\in H.$ By Proposition \ref{Orlicz_coro_BS_conjecture}, $\set{\overline{m_n}\,\|x_n\|\,y_n}$ is a Bessel sequence. Let $K=\sup\|x_n\|$. Since $\|x_n\|^p\leq K^{p-1}\|x_n\|$, it follows that $\set{\overline{m_n}\,\|x_n\|^p\, y_n}$ is a Bessel sequence.
\end{proof}

Using Theorem \ref{norm_Bes_orth_de} and Proposition \ref{Orlicz_coro_BS_conjecture}, we provide an example of Bessel-normalizable sequence for which the Stoeva-Balazs conjecture holds in the following.
\begin{example}
	Fix positive integers $L$ and $N$ and let $\set{e_n}$ be an orthonormal basis for $H$.
	For each $i\geq 0$, let $\set{x_n^i}_{n\Eq1}^L$ be an arbitrary sequence in $\clspan{(e_n)}_{n\Eq iN+1}^{i(N+1)}$ that does not contain any zero elements.
	Then by Theorem \ref{norm_Bes_orth_de}, we see that $\cup_i\set{x_n^i}_{n\Eq 1}^L$ is Bessel-normalizable. Assume that $\set{y^i_n}_{n=1}^L$'s are sequences in $H$ and $(m^i_n)_{n=1}^L$'s are sequences of scalars such that $\sum_i\sum_n m_n^i \,\ip{x}{y^i_n}\,x^i_n$ converges unconditionally for all $x\in H$. By Proposition \ref{Orlicz_coro_BS_conjecture}, the Stoeva-Balazs conjecture holds for $\sum_i\sum_n m_n^i \,\ip{x}{y^i_n}\,x^i_n$ and \text{\large$\cup_{i,\, n}$} $\bigset{\overline{m_n^i}\,\norm{x_n^i}\,y_n^i}$ is a Bessel sequence. \qeddef
\end{example}

\subsection{Frame-Normalizability of Iterative Systems Generated By Normal Operators}

We apply our results to partially answer the open question regarding the frame property of the iterative system $\bigset{\frac{A^nx}{\|A^n x\|}}_{x\in S,\, n\geq 0}$, where $A\colon H\rightarrow H$ is a normal operator and $S$ is a countable subset of $H$.

In the following, the notations $c_n \searrow L$ and $c_n \nearrow L$ will mean that the sequence of scalars $(c_n)$ monotonically decreases and monotonically increases to $L$ as $n\rightarrow \infty$, respectively.

We will need the following standard result.
\begin{lemma} 
	\label{fixed_pt_AandAstar}
	Let $A\colon H\rightarrow H$ be a bounded linear operator with $\norm{A}=1$. Assume that $x_0\in H$ is such that $Ax_0=x_0$. Then $A^*x_0=x_0$.
	\end{lemma}
\begin{proof}
	Without loss of generality, we may assume that $x_0\neq 0.$ Write $A^*x_0=ax_0+x_\perp$, where $a$ is a scalar and $x_\perp$ is the orthogonal projection of $x_0$ onto $\clspan{(x_0)}^\perp$. We compute that 
	$$\ip{x_0}{x_0}\Eq\ip{x_0}{Ax_0}\Eq \ip{A^*x_0}{x_0} \Eq a\, \ip{x_0}{x_0},$$
	which implies that $a=1$. Since $\norm{A^*x_0}\leq \norm{x_0}$, it follows that $x_\perp=0.$ Thus, $A^*x_0=x_0$.
	\end{proof}

\begin{theorem}
	\label{chara_normal_finite}
	Let $A\colon H\rightarrow H$ be a normal operator and let $S=\set{x_i}_{i=1}^N$ be a finite subset of $H$. Assume that one of the following conditions hold. 
	\begin{enumerate}
		\setlength\itemsep{0.3em}
		\item [\textup{(a)}] $\set{A^nx}_{x\in S,\, n\geq 0}$ is a frame.
		\item [\textup{(b)}] There exist $0<B\leq C$ such that $B\leq \norm{A^n x}\leq C$ for every $n\geq 0$ and $x\in S$.
		\item [\textup{(c)}] $\norm{A}=1$ and $A$ has a nonzero fixed point $w_0$ for which $\ip{w_0}{x_i}\neq  0$ for all $1\leq i\leq N.$
		\item [\textup{(d)}] $\norm{A}<1$ and $\set{P_M A^nx}_{n\geq0, x\in S}$ satisfies a lower frame condition in $M$ for some infinite-dimensional closed subspace $M$.
	\end{enumerate}
	Then $\bigset{A^nx}_{x\in S, \,n\geq 0}$ is not frame-normalizable. 
\end{theorem}
\begin{proof}
	(a) We mention that the first half of this proof closely follows the proof of \cite[Theorem 7]{AP17}. For the sake of completeness, we include the details here.
	Let $L$ and $K$ be frame bounds for $\set{A^nx}_{x\in S,\, n\geq 0}$. Take any $x_i\in S$, then for any $m\in \mathbb{N}$ we have $$\sum_{n=0}^\infty \sum_{x\in S}\,\bigabs{\bigip{(A^m)^*x_i}{A^n x}}^2\Eq \sum_{n=m}^\infty \sum_{x\in S}|\ip{x_i}{A^nx}|^2\Eq \sum_{x\in S}\sum_{n=m}^\infty |\ip{x_i}{A^nx}|^2  \Le  K\|x_i\|^2\,<\,\infty.$$ 
	Consequently,
	$\sum_{n=0}^\infty \sum_{x\in S}|\ip{(A^m)^*x_i}{A^n x}|^2$ converges to $0$ as $m\rightarrow \infty$.
	Since $$\sum_{n=0}^\infty \sum_{x\in S}|\ip{(A^m)^*x_i}{A^n x}|^2\Ge L\,\|(A^m)^*x_i\|^2,$$ we see that $\|(A^m)^*x_i\|\rightarrow 0$  as $m\rightarrow \infty$. Since $A$ is normal, it follows that $\|A^m x_i\|\rightarrow 0$.
	Becasue $S=\set{x_i}_{i=1}^N$ is finite, we can set $\set{y_n}=\set{x_1,\dots,x_N,Ax_1,\dots Ax_N,A^2x_1,\dots,A^2x_N,\dots}$ and obtain $y_n\rightarrow 0$ as $n\rightarrow \infty$. 
	By Corollary \ref{non_norma_frame_with_norm_conv_zero}, $\set{y_n}$ (and hence $\set{A^nx}_{x\in S,\, n\geq 0}$) is not frame-normalizable. (In fact, it is not even Bessel-normalizable.)
	
	\medskip
	
	(b)	Assume that $\Bigset{\frac{A^nx}{\|A^nx\|}}_{x\in S,\, n\geq 0}$ is a frame with frame bounds $L$ and $K$. The assumption consequently implies that $$\frac{L}{C}\,\norm{y}^2\Le \sum_{n=0}^\infty \sum_{x\in S}|\ip{y}{A^n x}|^2 \Le \frac{K}{B}\,\norm{y}^2$$ 
	Thus, $\set{A^nx}_{x\in S,\, n\geq 0}$ is a frame for $H$, which is a contradiction to (a).
	
	\medskip
	
	(c)	Since $\|A\|=1$,  it follows that $\|A^{n+1} x\|\leq \|A^n x\|$ for all $n\geq 0$ and all $x\in S$. Suppose that $\|A^nx\|$ does not decrease to $0$ for all $x\in S$. Then we have $$0<\min_{x_i\in S} \,\Bigparen{\lim\limits_{n\rightarrow \infty} \norm{A^nx_i}} \Le \norm{A^nx}\Le \max_{x_i \in S}\norm{x_i},$$ for all $n\geq 0$ and all $x\in S$. Thus, $\Bigset{\frac{A^nx}{\|A^nx\|}}_{x\in S,\, n\geq 0}$ is not a frame for $H$ by assumption (b). 
	
	Now assume there exists some $x_{i_0}\in S$ for which $\|A^nx_{i_0}\|\searrow 0$. By Lemma \ref{fixed_pt_AandAstar} and assumption (c), there exists some nonzero $w_0$ for which $(A^*)^n w_0=w_0$ for all $n\geq 0$. Let $K$ be an upper frame bound for $\Bigset{\frac{A^nx}{\|A^nx\|}}_{x\in S,\, n\geq 0}$. 
	Since $\ip{w_0}{x_{i_0}}\neq 0$, it follows that
	$$K\,\|w_0\|^2\Ge \sumli \Bigabs{\Bigip{(A^*)^nw_0}{\frac{x_{i_0}}{\|A^nx_{i_0}\|}}}^2\Eq |\ip{w_0}{x_{i_0}}|^2\sumli \frac{1}{\|A^nx_{i_0}\|^2}\Eq \infty.$$
	Thus, $\Bigset{\frac{A^nx}{\|A^nx\|}}_{x\in S,\, n\geq 0}$ is not a frame.
	
	\medskip
	
	(d) Since $\norm{A}<1$, it follows that $(\,\|A^nx\|\,)_n$ decreases to $0$ for all $x\in S$. Therefore, $\Bigset{\frac{A^nx}{\|A^nx\|}}_{x\in S,\, n\geq 0}$ is not frame for $H$ by Theorem \ref{non_norma_decrea_norm}.~(In fact, it is not even a Bessel sequence.)
\end{proof}

Next, we consider the case that $S$ contains only one element. 
\begin{lemma}
	\label{normal_incr_norm_lemma}
	Assume that $A\colon H\rightarrow H$ is a normal operator and fix $k_0\geq 0$.
	Then for any $n\geq 2$,  we have that $\displaystyle \|A^{k_0+n}x\|\Ge	\Bigparen{\frac{\|A^{k_0+1}x\|}{\|A^{k_0}x\|}}^{n-1}\|A^{k_0+1}x\|.$
	
	Therefore, if $x_0 \in H$ is such that $\|A^{k_0+1}x_0\|> \|A^{k_0}x_0\|$ for some $k_0\geq 0$, then $$\|A^{k_0+n}x_0\| \nearrow \infty\qquad \text{as $n\rightarrow\infty$.}$$ \ 
\end{lemma}
\begin{proof}
	We proceed by induction. Take any $x\in H$ with $\|x\|=1$. By the Cauchy-Bunyakovski-Schwarz inequality, we have that \begin{equation}\label{normal_op_ineq}
	\|Ax\|^2=|\ip{A^*Ax}{x}|\Le \|A^*Ax\|=\|A^2x\|. 
	\end{equation}
	By inserting $x=\frac{A^{k_0}x}{\|A^{k_0}x\|}$ into inequality (\ref{normal_op_ineq}), we see that for every $x\in H$ we have that  \begin{equation}\label{normal_op_ineq2}
	\|A^{k_0+2}x\|\geq	\frac{\|A^{k_0+1}x\|^2}{\|A^{k_0}x\|}\Eq \frac{\|A^{k_0+1}x\|}{\|A^{k_0}x\|}\|A^{k_0+1}x\|. 
	\end{equation}
	This establishes the base case $n=2$.

	Now assume that $\|A^{k_0+n}x\|\geq	\Bigparen{\frac{\|A^{k_0+1}x\|}{\|A^{k_0}x\|}}^{n-1}\|A^{k_0+1}x\|$ and $\frac{\|A^{k_0+n}x\|}{\|A^{k_0+n-1}x\|}\geq \frac{\|A^{k_0+1}x\|}{\|A^{k_0}x\|}$ for some $n\geq 2$. 
	By considering $x=A^{n-1}x$ in inequality (\ref{normal_op_ineq2}), we compute that 
	$$\|A^{k_0+n+1}x\|\Eq \|A^{k_0+2}(A^{n-1}x)\|\Ge \frac{\|A^{k_0+n}x\|}{\|A^{k_0+n-1}x\|}\|A^{k_0+n}x\|\Ge \biggparen{\frac{\|A^{k_0+1}x\|}{\|A^{k_0}x\|}}^{n}\|A^{k_0+1}x\|.$$
	The result then follows by induction.
\end{proof}
\begin{lemma}
	\label{norm_condition_normal_system_frame}
	Let $A\colon H\rightarrow H$ be a normal operator and fix $x_0\in H$. If $\Bigset{\frac{A^nx_0}{\|A^nx_0\|}}_{n\geq 0}$ is a frame for $H$, then there either exists some $k_0\in \N$ such that $\|A^{n+k_0}x_0\| \nearrow \infty$ or $\|A^nx_0\|\searrow 0$
	as $n\rightarrow\infty$.
\end{lemma}
\begin{proof}
	By Lemma \ref{normal_incr_norm_lemma}, $\|A^{n+k_0}x_0\|\nearrow \infty$ as $n\rightarrow \infty$ if there exists $k_0\in \N$ such that $\|A^{k_0+1}x_0\|\geq \|A^{k_0}x_0\|$. Assume that $\|A^{n+1}x_0\|\leq \|A^{n} x_0\|$ for every $n\geq 0$. If $\|A^nx_0\|$ does not decrease to zero, then we have $0< \lim\limits_{n\rightarrow \infty}\norm{A^nx_0}\leq \|A^{n} x_0\| \leq \norm{x_0}$ for all $n$, which is a contradiction to Theorem \ref{chara_normal_finite} (b).
\end{proof}

\begin{theorem}
	
	\label{norma_system_cant_frame_infi_proj}
	Let $A\colon H\rightarrow H$ be a normal operator and let $x$ be an element in $H$. Assume that one of the following conditions holds.
	\begin{enumerate}
		\setlength\itemsep{0.3em}
		\item [\textup{(a)}] There exists an infinite-dimensional closed subspace $M$ such that $\set{P_MA^nx}_{n\geq 0}$ is a frame for $M$.
		\item [\textup{(b)}] $\|A\|\leq 1$ and there exists an infinite-dimensional closed subspace $M$ such that $\set{P_MA^nx}_{n\geq 0}$ satisfies a lower frame condition in $M$.
		\item [\textup{(c)}] There exists a finite subset $J$ of $\N$ such that $\set{A^nx}_{n\notin J}$ is a frame for its closed span.
	\end{enumerate}
	Then $\Bigset{\frac{A^nx}{\|A^nx\|}}_{n\geq 0}$ is not a frame for $H$.
\end{theorem}
\begin{proof}
	(a) and (b)
	Let $M$ be any infinite-dimensional closed subspace of $H$.
	Suppose that $\Bigset{\frac{A^nx}{\|A^nx\|}}_{n\geq 0}$ is a frame. By Lemma \ref{norm_condition_normal_system_frame}, we either have $\|A^nx_0\| \searrow 0$ or $\|A^nx_0\| \nearrow \infty$. If $\|A^nx_0\| \searrow  0$ (in particular, this is the case if $\|A\|\leq 1$), then $\set{P_MA^nx}_{n\geq 0}$ does not satisfy a lower frame condition by Theorem \ref{non_norma_decrea_norm}.  On the other hand, $\set{P_MA^nx}_{n\geq 0}$ cannot be a Bessel sequence by Theorem \ref{non_norma_increa_norm} if $\|A^n x_0\| \nearrow \infty$. Therefore, $\set{P_MA^nx}_{n\geq 0}$ is not a frame for $M$.
	
	\medskip
	(c)	Let $M=\overline{\text{span}}\set{A^nx}_{n \notin J}$. 
	If $M$ is finite-dimensional, then $\Bigset{\frac{A^nx}{\|A^nx\|}}_{n\geq 0}$ is not a frame by Proposition \ref{fini_dim_char}. Thus, we may assume that $M$ is infinite-dimensional. Suppose that $\set{A^nx}_{n\notin J}$ is a frame sequence. The union of a frame for $M$ and any finite subset of $M$ is still a frame for $M$. Consequently, $\set{P_M A^n x}_{n\geq 0} = \set{P_M A^n x}_{n\in J} \cup \set{A^nx}_{n\notin J} $ is a frame for $M$, which contradicts statement (a).
\end{proof}

\subsection{Frame-Normalizability of Iterative Systems Generated By Compact Operators}
Let $A\colon H\rightarrow H$ be a compact operator and let $S$ be a finite subset of $H$.
It was proved by Christensen et al.\ in \cite{CHR18} that $\set{A^nx}_{x\in S,\, n\geq 0}$ cannot be a frame for $H$. We close this paper by giving some sufficient conditions such that $\Bigset{\frac{A^nx}{\|A^nx\|}}_{x\in S,\, n\geq 0}$ cannot be a Bessel sequence in $H$.
\begin{proposition}
	Let $A\colon H\rightarrow H$ be a compact operator and let $S\subseteq H$ be a finite subset of $H$. Assume that one of the following conditions holds.
	\begin{enumerate}
			\setlength\itemsep{0.3em}
		\item [\textup{(a)}] There exists an infinite dimensional closed subspace $M$ of $H$ such that $\set{P_MA^nx}_{x\in S,\, n\geq 0}$ satisfies a lower frame condition in $M$.
			\item [\textup{(b)}] There exists a constant $C>0$ such that $ \|A^n x\|\Ge C$ for every $n\geq 0$ and $x\in S$.
		   \item [\textup{(c)}] $\norm{A}=1$ and $A$ has a nonzero fixed point $w_0$ for which $\ip{w_0}{x}\neq  0$ for all $x\in S.$
	\end{enumerate}
Then $\set{A^nx}_{x\in S,\, n\geq 0}$ is not Bessel-normalizable.
\end{proposition}    
\begin{proof}
	Fix $x_0\in S$. We show that if $\Bigset{\frac{A^nx}{\|A^nx\|}}_{x\in S,\, n\geq 0}$ is a Bessel sequence for $H$, then $\norm{A^n x_0 }$ converges to $0$ as $n\rightarrow\infty$. Let $B$ be an upper frame bound for $\Bigset{\frac{A^nx}{\|A^nx\|}}_{n\geq 0}$. Since $\sum\, \bigabs{\ip{y}{\frac{A^n x}{\norm{A^n x}}}}^2\leq B\,\|y\|^2$ for every $y\in  H$, we have that $\bigset{\frac{A^n x}{\norm{A^n x}}}$ converges weakly to $0$. By the fact that compact operators send weakly convergent sequences to norm-convergent sequences, we see that $\bignorm{\frac{A^{n+1} x}{A^n x}}$ converges to $0$. Therefore, there exists $n_0\in \N$ such that $\norm{A^{n+1}x}\Le \frac{\norm{A^nx}}{2}$ for any $n\geq n_0$. Thus, $\norm{A^{n_0+k}x}\Le 2^{-k} \norm{A^{n_0}x}$ for every $k\in \N$ and hence $\,\norm{A^n x} \rightarrow 0$. 
\medskip
	
	(a) Arguing similarly to Theorem \ref{chara_normal_finite}(a), the result then follows by Theorem \ref{non_norma_decrea_norm}.
	\medskip
	
	(b) Since $\norm{A^n x}$ converges to 0 if $\Bigset{\frac{A^nx}{\|A^nx\|}}_{x\in S,\, n\geq 0}$ is a Bessel sequence, it follows that $\set{A^nx}_{x\in S,\, n\geq 0}$ is not Bessel-normalizable if it is norm-bounded below.
	\medskip
	
	(c) This follows by the same argument we used for Theorem \ref{chara_normal_finite}(c).	
	\end{proof}

\section*{Acknowledgments}
We thank Victor Bailey for bringing reference \cite{ACMCP16} into our attention, which is an inspiration for some of our future work. We would like to express gratitude to Christopher Heil for fruitful discussions and helpful comments on the draft.

\end{document}